\documentclass[12pt]{article}
\usepackage{newinutile}
\usepackage{amsmath,amssymb}
\usepackage{epsfig,graphics,color,calc}
\usepackage{cite}

\newcommand{\qed}{\hfill $\Box$\bigskip}
\newcommand{\rr}[1]{{\normalfont\textrm{#1}}}
\newcommand{\bb}[1]{{\mathbb{#1}}}
\newcommand{\cR}[1]{{\mathcal{R}}}

\newcommand{\resist}{\widetilde\Gamma}

\newcommand{\gb}{\beta}
\newcommand{\gep}{\varepsilon}

\newlength{\pecettawidth}
\begin{document}

\title{A comparison between different cycle decompositions for
Metropolis dynamics}

\maketitle

\author{Emilio N.M.\ Cirillo}
\affiliation{Dipartimento di Scienze di Base e Applicate per l'Ingegneria
(SBAI), Sapienza Universit\`a di Roma,
via A.\ Scarpa 16, I--00161, Roma, Italy.}
\email{emilio.cirillo@uniroma1.it}
\thanks{ENMC acknowledges the Mathematics Department of the
Utrecht University for the kind hospitality and for financial support.}

\author{Francesca R.\ Nardi}
\affiliation{Department of Mathematics and Computer Science,
             Eindhoven University of Technology,
             P.O.\ Box 513, 5600 MB Eindhoven, The Netherlands.}
\affiliation{Eurandom, P.O.\ Box 513, 5600 MB, Eindhoven, The Netherlands.}
\email{F.R.Nardi@tue.nl}

\author{Julien Sohier}
\affiliation{Department of Mathematics and Computer Science,
             Eindhoven University of Technology,
             P.O.\ Box 513, 5600 MB Eindhoven, The Netherlands.}
\affiliation{Eurandom, P.O.\ Box 513, 5600 MB, Eindhoven, The Netherlands.}
\email{j.sohier@tue.nl}


\begin{abstract}
In the last decades the problem of metastability has been attacked on
rigorous grounds via many different approaches and techniques which are briefly reviewed in this paper.
It is then useful to understand connections between different
point of views. In view of this
we consider irreducible, aperiodic and
reversible Markov chains with exponentially
small transition probabilities in the framework of Metropolis dynamics.
We compare two different cycle decompositions and prove their equivalence.
\end{abstract}




\section{Introduction}
\par\noindent
Cycle decomposition is a very useful tool to estimate first hitting
times for stochastic processes. In this note we compare, and prove
the equivalence, between two different approaches that will be
respectively called
\textit{path} and \textit{graph cycle} decompositions.
These results are utterly important in the generic study of the metastability phenomenon.

\subsection{Metastability}
\par\noindent
The phenomenon of metastability is defined by the following scenario:
(i) a system is ``trapped'' for an abnormally long time in a state ---
the \textit{metastable state} --- different from the eventual equilibrium
state consistent with the thermodynamical parameters.
(ii) In the metastable state the system behaves as if it were in
regular equilibrium.
(iii) Subsequently, the system undergoes a \emph{sudden} transition at a
\emph{random time} from the metastable to the stable state.
The mathematical study of this phenomenon has been a standing issue since
the foundation of Statistical Mechanics, but only in the 80's
rigorous mathematical approaches started to be developed and, due to
the great interest of the subject, they then proliferated to a multitude of
different approaches.
These mathematical approaches, however, are not completely
equivalent as they rely on
different definitions of metastable states and thus involve
different properties of hitting and escape times.
The situation is particularly delicate for evolutions of infinite--volume
systems and for irreversible systems.
The proliferation of definitions and hypotheses on metastable behavior
arises from the diversity of the physical situations in which the phenomenon
appears. This diversity results in contrasting demands on the appropriate
mathematical theory. The main issues confronted by the mathematical
treatment of metastability can be grouped into three categories.

\textit{Conservative vs non--conservative dynamics.\/}
This dichotomy applies to dynamics for statistical mechanical models of
fluids or magnets.
Non--conservative dynamics are those that do not conserve the total number
of particles or the total magnetization. They include Glauber (that is,
single spin--flip) dynamics (used to model metastable ferromagnets)
and many probabilistic cellular automata (that is, parallel dynamics).
In contrast, conservative dynamics are suitable to the study of supersaturated
gases. Its study poses {\it enormous challenges}
because particle or magnetization conservation introduces non--local effects.

\textit{Finite vs infinite configuration space.\/}
Two extreme types of metastability studies can be distinguished.
In the \emph{finite--space} case, the configuration space remains fixed
(or bounded) while the drift towards (meta)stable configurations is
increased (e.g., temperature goes to zero). In the \emph{infinite--space}
setup, the size of the configuration space grows in an unbounded fashion,
while drifts are kept approximately constant.
In many instances both parameters (drift and size of the space)
change simultaneously, but usually these changes are coupled so one of the
effects is dominant. Mathematically, the distinction stems from the
possibility of ``entropic'' effects in the infinite--volume case that
changes the scale at which distribution laws must be controlled.
The iconic case is the thermodynamic limit of spin or gas models,
see for instance \cite{SS,dHOS,GHNOS09,GN}.
In these models, exit from metastability requires \emph{nucleation},
that is the formation of a critical droplet.
The probability of such an event in a large volume must include the
``entropic contribution'' due to the fact that the nucleation can take
place anywhere in the volume.

\textit{Parallel dynamics and cost functions.\/}
Following the matrix imposed by Metropolis and Glauber dynamics,
stochastic transition rates are written as exponentials of
\emph{cost functions}. For reversible single spin--flip dynamics these costs
are determined by the difference of energy between the two states involved
in the move. This is not so for parallel dynamics (for instance for
Probabilistic Cellular Automata \cite{BCLS,CN,CNS}) in which, at each step,
all spins are independently tested. In such evolutions, costs are a possibly
complicated function of the different patterns of spin flips connecting
the relevant states. In these cases it is therefore necessary to
dissociate energy profiles from \emph{energy barriers}.
The former are associated to invariant measures and determine the set of
stable and metastable configurations. The latter are associated to
transition rates and determine drifts and exit times.

\subsection{Different approaches to metastability}
\par\noindent
Early approaches to metastability were based on the computation of
expected values with respect to restricted equilibrium
states \cite{PL}.  
This paradigm is still influential in physics, see e.g. \cite{LL}.
The influence of Probability on Statistical Mechanics led to
an alternative pattern of thoughts based on stochastic dynamics and focused
on their spectral properties and on the behavior of their typical trajectories.
This point of view has given rise to different theoretical constructions that
can be classified, roughly, in three major groups.

\textit{(I) Classical approach: Hitting times of Markov chains.\/}
The escape time from metastability is determined, by the \textit{visiting}
or \textit{hitting} time to a set of states of very small (invariant)
measure, when most of the measure is carried by a different, somehow reduced,
set of states (see, e.g. \cite{Kei79}).
Similar problems were confronted in reliability theory where typical states
were called \emph{good} and those concerned by the hitting times were
called \emph{bad}. The exponential character of good--to--bad transitions
is due to the existence of two different
time scales: long times are needed to go from good to bad states, while the
return to good states from anywhere (except, perhaps, the bad states)
is much shorter. As a result, a system in a good state can reach the bad ones
only through a large fluctuation bringing it all the way to the bad state.
Indeed,
any ``intermediate" fluctuation would be followed by an unavoidable return
to the good states, where, by Markovianity, the process would start afresh
independently of previous attempts.
Hence,
the escape time is made of a large number of independent returns
to the good states followed by a final successful excursion to badness
that must happen without hesitation, in a much shorter time.
The exit time is, therefore, a geometric random variable with extremely
small success probability. In the limit exponentiality is found.


\textit{(II) Large deviations of trajectory distributions.\/}
Freidlin and Wentzel \cite{FW} were the first to use the large
deviation machinery to study the problem of exit times from an
attractive domain. Their theory applies to stochastic differential
equations with a deterministic driving gradient force and a small
Brownian stochastic term. The deterministic part of  the dynamics
is responsible for the fast return to ``good states'' while the
stochastic contribution provides the escape mechanism.
The smallness of this last contribution leads to very long time scales
for the visit to ``bad'' states.
Typical trajectories are described using a graphical method
built out of two basic ingredients:
\emph{cycles} (associated to metastable pieces of trajectories)
and \emph{exit tubes} (describing typical escape trajectories).

The Freidlin and Wentzel theory evolved into two related schools that
we shall call the
\emph{graphical} \cite{CaCe} and the \emph{pathwise} \cite{CGOV} {approach}.
The former relies on a refinement of Freidlin and Wentzel's graphical methods
allowing for a detailed study of exit paths via a decomposition
into cycles and saddle points traversed from cycle to cycle.
The exit time also decomposes into the time spent at each point of the exit
path. This graphical approach has been applied to reversible Metropolis
dynamics and to simulated annealing in \cite{Ca,CaCe,CaTr,Trouve1,Trouve,Trouverl}.

The pathwise approach, on the other hand, proposed in \cite{CGOV},
was introduced as an adaptation of the ideas of Freidlin
and Wentzel to Metropolis--like dynamics, with all notions and properties
expressed in terms of an energy profile associated to the invariant measure.
This provides a clearer and physically appealing picture.
In particular, absolute energy minima identify stable states and ``deep" local
energy minima lead to metastability.
The two time scales within each energy well correspond, respectively, to fast
``downhill'' and infrequent ``uphill'' trajectories.
In the limit of very steep wells (temperature tending to zero), the theory
yields rather precise information on:
(i)
the \emph{transition time}, i.e., the time needed to arrive to
the stable equilibrium, which is determined by the height of the largest
energy barrier separating metastable from stable states.
(ii)
The \emph{typical exit tube}, i.e., the sequences of configurations
along which the formation of the stable phase takes place.
This is the physically relevant mechanism that, in gas or spin systems,
is mediated by the appearance of a {\it critical droplet} after which the
system quickly relaxes to equilibrium.

The full power of this method was first exploited in \cite{OS1} and
comprehensively reviewed in \cite{OV}.
It has been extended to non--reversible Markov processes in \cite{OS2}
(though irreversibility brings back to the graphical
approach). The approach was further
simplified in \cite{MNOS} where transition times are determined on the basis
of a ranking of \emph{stability levels}, without requiring detailed
knowledge of typical trajectories.

\textit{(III) Potential--theoretical approach:
spectral properties of Markov transition matrices.\/}
In the early eighties Aldous and Brown \cite{AB1} proposed a
new approach to the hitting--time theory based on spectral properties of
the transition matrix and the use of the Dirichlet form.
This approach has the advantage of leading to quite precise error bounds for
the exponential approximation. The current version of this strategy is the
potential--theoretic approach developed
in \cite{BEGK01} 
(see \cite{BovMetaPTA} for reviews).
Besides exponential laws, this method gives more precise estimates
of the expected value of the
transition time, including a prefactor that cannot be
found with alternative approaches. The determination of this prefactor,
however, requires the knowledge of the critical droplet and neighboring
configurations; information that has to be imported from more detailed
pathwise studies.
In \cite{BL1} another use of spectral and potential theoretical
techniques is proposed in which only visits to well bottoms are registered.
Upon time rescaling, a continuous--time Markov process is obtained whose
transition rates encode the information on transition times. 
(See also \cite{BLM,BG} for recent development).

\subsection{Application overview}
\par\noindent
We outline briefly some applications of the theories
described above.
The aim is not to be exhaustive, but rather to list references useful in
relation to the definitions and comparisons to follow.

The general theory \cite{OS1,OS2,Scop94,CaCe,OV,MNOS,Ca} and metastability studies
in the nineties (see
\cite{CGOV,NeS,NevSchbehavdrop,CN12,KO2,CO,AC,NO,CL} for the pathwise approach
and \cite{CaTr,Trouve1,Trouve,Trouverl} for the graphical approach)
refer to single--spin flip dynamics (Metropolis, Glauber)
of Ising--like models (including mean--field versions) in finite volume
and at low temperature.
Studies within the general potential theoretical approach
refer both to a general Markov--process point of view
\cite{BEGK02,BEGK04} and applications to
mean--field and Ising model \cite{BEGK01,BM}.

The study of metastability for conservative dynamics started a decade after
and initially involved nearest--neighbor lattice gases at low temperature
and density inside finite boxes with open boundary conditions
mimicking infinite reservoirs fixing particle density.
Relevant references are \cite{dHOS,GOS,HNOS,NOS} for the pathwise approach
and \cite{BHN,HNT1,HNT2,HNT3} for the potential theoretic approach.

Models with parallel dynamics
were studied first from a numerical point of view in \cite{BCLS} and then
rigorously in \cite{CN,CNS,CNP,CNS2,CNS3} (pathwise approach) and
\cite{NS1} (potential theoretic approach).

The more involved infinite--volume limit, at low temperature or vanishing
magnetic field, was first studied via large deviations techniques in
\cite{DS2,S2,MO1,MO2,SS,CeMa} and potential
theoretically in \cite{BHS}. These references dealt with Ising and
Blume--Capel models under Glauber dynamics.
The Ising lattice gas model subjected  to Kawasaki dynamics was studied
in \cite{dHOS,GHNOS09,GN} and \cite{BHS} (potential theoretic approach)
in the limit of temperature and volume growing exponentially fast to infinity.

\subsection{Aim of the paper}
\par\noindent
As it has been explained above, due to the great theoretical and
applicative interest of metastability, different mathematical
theories of metastability have been developed in the past years.
It is interesting to understand the mutual connections
in order to apply results proven in one framework to systems
that are naturally approached in a different one.

We have also remarked that within each theory different
flavors have appeared. In particular,
in the framework of the ``large deviation" point of view,
two different approaches to cycle decomposition, the \textit{graphical}
and the \textit{path} one, have been developed.
The former has been introduced and applied in a very general setup
while the latter as been first introduced for Metropolis dynamics
and then extended to models with cost functions such as
Probabilistic Cellular Automata \cite{CN,CNS}. Due to its generality,
the results on hitting times proven with the graphical approach
are written in terms of complicated expression whose physical
meaning is sometimes difficult to be caught. On the other hand, in the
framework of the path approach everything is clearly written in terms
of differences of energies, but the theory applies only to Metropolis--like
systems.

Our opinion is that it is interesting and it can also be very
fruitful to understand the connections between these two approaches
in order to mix the strength of the former with the simplicity of the latter.
As a first step in this direction we prove the equivalence of the
two approaches in the case of the Metropolis dynamics, whose definition is
now recalled.

\subsection{Setup}
\label{s:setup}
\par\noindent
We consider a finite state space $S$ equipped with a  function
$H: S \rightarrow \bb{R}$. Sometimes for a point $x \in S$, we will write the \textit{energy} of $x$ for the value $H(x)$.

  We assume that $S^{2}$ is equipped with a
connectivity function $q: S^{2} \rightarrow [0;1]$, which satisfies the
following conditions:
(i) for any $x \in S, \sum_{y \neq x} q(x,y) \leq 1$;
(ii) for any $(x,y) \in S^{2},  q(x,y) = q(y,x)$;
(iii) for any  $(x,y) \in S^{2}$, there
exists $n \in \bb{N}$ and $x_{1}, \ldots, x_{n}$
such that $x_{0} = x, x_{n} = y$ and $q(x_{i},x_{i+1}) >0$
for $i = 0,\ldots, n-1$.

For $\gb>0$, we then define the
\textit{Metropolis Markov chain} $X$ as being the Markov chain with
transitions given by $p_{\gb}$ where the kernel $p_{\gb}$ satisfies
\begin{equation*}
 p_{\gb}(x,y) = q(x,y) \exp\left( -\gb \left( H(y) - H(x) \right)^{+} \right)
\end{equation*}
if $x \neq y$ and $p_{\gb}(x,x) = 1 - \sum_{y \neq x} p_{\gb}(x,y)$.

The chain started at $x\in S$ will be denoted by $x_0=x,x_1,\dots,x_t,\dots$
and
the associated probability will be denoted by
$\mathbb{P}_x$.
The main notion the paper will deal upon is that of \textit{hitting} time
\begin{equation}
\label{hitting}
\tau_G:=\inf\{t,\,x_t\in G\}
\end{equation}
to a set $G\subset S$.

A particularly famous example of Metropolis dynamics is given by the standard
Ising model under Glauber dynamic.

The purpose of this paper is to discuss the first exit problem of a general
set $G$ for the dynamics defined above for large $\gb$. It turns out that
the most relevant case corresponds to $G$ given by a cycle, namely a set whose
internal points, for large $\gb$, are typically visited many times by our
process before exiting.

The paper is organized as follows:
in Section~\ref{s:path} and Section~\ref{s:FWC}
we recall the definition of cycles and the main results on hitting
times respectively in the framework of the path and the graph approach.
In Section~\ref{s:equiv} we explore the connections between the
two approaches and state their equivalence.

\section{Path cycles}
\label{s:path}
\par\noindent
We briefly review the path approach to cycle decomposition in \cite{OS1}.
In particular we discuss few properties that will be useful
in Section~\ref{s:equiv}.
A path $\omega$ is a sequence $(\omega_{1}, \ldots, \omega_{n})$
of communicating states, that is to say $q(\omega_{i}, \omega_{i+1}) > 0$
for any $i = 0, \ldots, n-1$.
We write $\omega: x \rightarrow y$ to denote a path joining $x$ to $y$.

We say that a subset $G$ of $S$ is \textit{connected} if for any $x,x' \in G$,
there exists a path $\omega: x \rightarrow x'$ such that $\omega$ is
entirely contained in $G$.
Two not empty subsets $G,G'\subset S$ are \textit{connected} whenever there exists
$x\in G$ and $x'\in G'$ such that $q(x,x')>0$.

If $G$ is a subset of $S$ on which $H$ is constant, we will write
$H(G)$ for the value of $H$ on $G$.

\begin{definition}
Let $G \subset S$, we define the exterior boundary $\partial G$ of $G$
and the ground $F(G)$ of $G$ respectively as
\begin{displaymath}
\partial G := \{ y \in S\setminus G, \exists x \in G, q(x,y) > 0 \} \\
\;\;\;\textrm{ and }\;\;\;
F(G) := \{ x \in G, H(x) = \min_{G} H \}
\end{displaymath}
\end{definition}

For a subset $G \subset S$ and $x \in S$, we say that $x$ is a \textit{neighbor} of $G$ if $x \in \partial G$.

\begin{definition}
\label{defO}
The set $A \subset S$ is a \textit{ non--trivial path cycle} of $S$ if and only if
 it is a connected subset of $S$ verifying
\begin{equation}
\label{defC}
\max_{A} H
<
\min_{\partial A} H
\end{equation}
We say that a subset $A \subset S $ is a cycle if and only if $A$ is a singleton or $A$ is a non--trivial path cycle.
\end{definition}

   In other words, a singleton is a trivial path cycle if and only if it is not a local minimum of $H$.

\begin{lemma}
\label{lvlset}
(Proposition 6.7 \cite{OV})
Given a state $x \in S$ and a real number $c\ge H(x)$,
the set of all points connected to $x$ by paths whose points
have energy smaller or equal to $c$ is
either the trivial path cycle $\{x\}$ or
is a non--trivial path cycle containing $x$.
\end{lemma}


\noindent
In the particular case $c = H(x)$ we denote by $U_{\leq x}$ the path cycle
whose existence is ensured by the above lemma. In words,
$U_x$ is the path cycle made of all the points in $S$ connected to $x$ via a
path whose points are at energy smaller or equal to $H(x)$.

\begin{lemma}
\label{interincl}
(Proposition 6.8 \cite{OV})
Let $A_1,A_2\subset S$ be two path cycles such that
$A_1 \cap A_2 \neq \emptyset$.
Then either $A_1 \subset A_2$ or $A_2 \subset A_1$.
\end{lemma}

\noindent
As a trivial consequence of these two lemmas, note that given a cycle
$A$ and $x \in A$, one has the inclusion $U_{\leq x} \subset A$.
On the other hand, as soon as $x \in A$ satisfies $H(x) = \max_{A} H$,
one gets the equality $U_{\leq x} = A$.

We remark, now, the following interesting property: two non--trivial
disjoint cycles cannot be connected. More precisely we state the following
lemma.

\begin{lemma}
\label{path}
Let $A_1,A_2\subset S$ be two path cycles. If $A_1$ and $A_2$ are
connected, then either $|A_1| = 1$ or $|A_2|=1$.
\end{lemma}

\par\noindent
\textit{Proof.\/}
By contradiction, assume that both $|A_1|>1$ and $|A_2|>1$
(and hence that they both satisfy equation \eqref{defC}).

Since the two path cycles
$A_1$ and $A_2$ are connected,
there exists $x_1\in A_1\cap\partial A_2$ and, thus,
it follows that
\begin{displaymath}
\max_{A_1} H \ge H(x_1) \ge \min_{\partial A_2} H > \max_{A_2} H
\end{displaymath}
where in the last bound we used \eqref{defC}).
But similarly, there exists $x_2\in A_2\cap\partial A_1$, which
implies
\begin{displaymath}
\max_{A_2} H \ge H(x_2) \ge \min_{\partial A_1} H > \max_{A_1} H
\end{displaymath}
which contradicts the above inequality.
\qed

\begin{definition}
\label{defdep}
Given a non--trivial path cycle $A$, we let
the \textit{depth} $\Gamma(A)$
and the \textit{resistance height} $\resist(A)$ of the cycle be respectively
\begin{equation*}
\Gamma(A): = \min_{\partial A} H - \min_A H
=H(F(\partial A))-H(F(A))
\end{equation*}
and
\begin{equation*}
\resist(A):=
\max_{A} H - \min_{A} H
=
\max_{A} H - H(F( A))
\end{equation*}
\end{definition}

\par\noindent
The following result describes the way the Markov chain exits a cycle
$A$ at very low temperature. As a matter of fact, it is known in a
more general, not reversible setup satisfying suitable hypotheses
called the Freidlin--Wentzell conditions.

\begin{theorem}
\label{mainOV}
(Theorem 6.23 in \cite{OV})
Given a non--trivial cycle $A$ and $\varepsilon > 0$,
for any $x,x' \in A$, the following properties hold in the asymptotic
$\beta \to \infty$:
\begin{equation}
\label{iptOV}
\bb{P}_{x}\left[\exp\{\beta(\Gamma(A)+\varepsilon)\}
                >
                \tau_{\partial A}
                >
                \exp \{\beta(\Gamma(A)-\varepsilon)\}\right]
        = 1 - o(1)
\end{equation}
and
\begin{equation}\label{iiptOV}
\bb{P}_{x}\big[\tau_{x'}<\tau_{\partial A},\;\;
                \tau_{x'}<\exp\{\beta(\resist(A)+\gep)\} \big]
          = 1 - o(1)
\end{equation}
\end{theorem}

\par\noindent
Roughly speaking, equation \eqref{iptOV} states that, starting from any
point of the cycle $A$, the exit time from $A$ is of order
$\exp\{\beta \Gamma(A)\}$ in the large $\gb$ asymptotic.
On the other hand, starting from any point in $A$, equation \eqref{iptOV}
says that, before exiting $A$, the Markov chain visits all the configurations
in $A$ within a time of order $\exp\{\beta \resist(A)\}$.

\section{Graph cycles}
\label{s:FWC}
\par\noindent
The construction of graph cycles due to Freidlin Wentzell \cite{FW}
is performed recursively.
Here, we recall this construction
following \cite[Part~2]{Trouverl}. In Section~\ref{Esemp} we 
discuss an example.

\subsection{Construction of graph cycles.}\label{reccons}
\par\noindent
Before starting with the recursive construction we need to recall
some general definitions: given a not empty set $M$ we denote by
$\mathcal{P}(M)$ the collection of all the subsets of $M$. Moreover,
a function $f:M\times M\to\mathbb{R}^{+} \cup \{ \infty \}$, namely, a function associating
each pair of elements
of $M$ with a (not necessarily finite) not negative real number will be
called a \textit{cost function}
on $M$. A \textit{path} of elements of $M$ is an element
$(m_1,\dots,m_n)\in M^n$ for some $n$ positive integer.
We shall misuse the notation by
also writing
\begin{displaymath}
f(m)=\sum_{i=1}^{n-1} f(m_i,m_{i+1})
\end{displaymath}
for any path $m=(m_1,\dots,m_n)\in M^n$.

  The following recursive construction can be read together with the example developped in Section \ref{Esemp}.

Recalling that the setup\footnote{It is important to remark
that graph cycles \cite{FW} are usually introduced in a more general
setup. For Markov chain it is usually assumed the so called
Freidlin--Wentzel assumption, namely,
there exists $\kappa>1$ such that
\begin{displaymath}
(1/\kappa) q(x,y) \exp\{-\beta V(x,y)\}
\leq p_{\gb}(x,y)
\leq \kappa q(x,y) \exp\{-\beta V(x,y)\}
\end{displaymath}
Note that the Metropolis dynamics is just a particular case.
}
is the one introduced in Section~\ref{s:setup}, we
define the \textit{zero--order set of graph cycles}
$E^{0} := \{ \{i\}, i \in S\}$ and 
the associated \textit{zero order cost
function} $V^{0}(\{i\},\{j\}) := (H(j) - H(i))^{+}$ if $i$ and $j$ 
are connected and 
$V^{0}(\{i\},\{j\}) := \infty$ otherwise.

Assume, then, that the $k$--order set of graph cycles
$E^{k} \subset \mathcal{P}(S)$ is constructed and equipped with
the $k$--order cost function $V^{k}$.
To implement the recursion, we proceed in five steps:
\begin{enumerate}
\item
for $A \in E^{k}$, let
 \begin{equation}\label{He}
 H_{\rr{e}}^{k}(A) := \min\{ V^{k}(A,A'), A' \in E^{k} \}.
 \end{equation}
Moreover, we define the renormalized cost function $V^{k}_{*}$ on $E^{k}$ by
setting
\begin{equation}\label{defVs}
V_{*}^{k}(A,B) := V^{k}(A,B) - H_\rr{e}^{k}(A)
\end{equation}
for all $A,B\in E^{k}$.
\item
For $A,A'$ elements of $E^{k}$, define the $\xrightarrow{k}$
relation by $A \xrightarrow{k} A'$ if and only if there exists a path
$\omega$ of elements of $E^{k}$
starting from $A$ and ending in $A'$ such that the cost of
$\omega$ with respect to $V_{*}^{k}$ is zero,
in short $V_{*}^{k}(\omega) = 0$.
\item
For $A,A'$ elements of $E^{k}$, define the relation of equivalence
$\mathcal{R}_{k}$ by $A \hspace{2pt} \mathcal{R}_{k} \hspace{2pt} A'$
if and only if both $A \xrightarrow{k} A'$ and $A' \xrightarrow{k} A$.
Then we stick together all the distinct classes of equivalence and
define the set
$D^{k+1} := \{ \bigcup_{A':A \mathcal{R}_{k} A'} A', A \in E^{k} \}.$
\item
On the set $D^{k+1}$, define the (partial) order
${\geq}^{k+1}$ by $A \geq^{k+1} A'$ if and only if there exist
$B,B' \in E^{k}, B \subset A, B' \subset A'$ such that
$B\xrightarrow{k}B'$.
Then we introduce $D_{^\star}^{k+1}$ the set of minimal
elements for the order relation $\geq^{k+1}$.
\item
Define $E^{k+1}$ as being the union of the set
$D_{^\star}^{k+1}$ and of the elements of $E^{k}$
which are not subsets of $D_{^\star}^{k+1}$, namely
\begin{equation}\label{En}
E^{k+1} := D_{^\star}^{k+1} \bigcup \{ A, A \in E^{k},
\exists B \in D^{k+1} \setminus D_{^\star}^{k+1}, A \subset B \}.
\end{equation}
For $A \in E^{k+1}$, define
 \begin{equation}\label{Hm}
 H_\rr{m}^{k+1}(A) := \max \{ H_\rr{e}^{k}(A'), A' \in E^{k}, A' \subset A \}.
 \end{equation}

\item
Finally define the cost function $V^{k+1}$ on $E^{k+1}$ by
\begin{equation}\label{defVn}
V^{k+1}(A,A') := H_\rr{m}^{k+1}(A) +
\min \{ V_{\star}^{k}(B;B'), B,B' \in E^{k}, B \subset A, B' \subset A' \}.
\end{equation}
\end{enumerate}


 The construction continues until $E^{k} = \{S\}$. As noted in \cite{Trouverl}, the recursive procedure described here is not stationary until iteration  
$n_{S}$, where $n_{S}$ is the first iteration such that $E^{n_{S}} = S$.

We remark that for any $k \geq 0$, $E^{k}$ is a partition of $S$ and more precisely the procedure gives a hierarchical decomposition of the state space
as a tree starting from the singletons and ending with the whole space.

We define the set of \textit{graph--cycles}
$\mathcal{C} := \bigcup_{k \geq 0} E^{k}$ and call
any
element of $\mathcal{C}$ a \textit{graph--cycle}.

\begin{definition}
Let $A \in \mathcal{C}$.
We set
   \begin{align*}
    H_\rr{e}(A)
    := \sup_{k \geq 0} \{ H_\rr{e}^{k}(A) \}
    \end{align*}
if $A\neq S$
     and $H_\rr{e}(S) = \infty$.
\end{definition}

The above definition is based on the following remark (see \cite{Trouverl}):
for $A \in \mathcal{C}$, it is easy to see that,
whenever $A \in E^{k} \cap E^{k+1}$, one has
\begin{equation}\label{drop}
H^{k+1}_\rr{m}(A) = H_\rr{e}^{k+1}(A) = H_\rr{e}^{k}(A)
\end{equation}
and thus $H_\rr{e}^{k}(A) = H_\rr{e}(A)$ as soon as $A \in E^{k}$.

\begin{definition} Let $A \in \mathcal{C}$ such that $A \neq S$ and
$|A|>1$. We set
  \begin{itemize}
 \item[--]
  $\mathcal{C}_{A}^{\star} := \{ B \in \mathcal{C}, B \subset A, B \neq A \}$;
  \item[--]
 ${\displaystyle
   \mathcal{M}_{\star}(A) := \{ B \in \mathcal{C}, B \hspace{2pt}
   \text{is a maximal element in} \hspace{4 pt} \mathcal{C}_{A}^{\star} \}
  }$
  (maximal proper partition of $A$);
 \item[--]
 ${\displaystyle H_\rr{m}(A) := \sup \{ H_\rr{e}(B),
                           B \in \mathcal{C}_{A}^{\star} \} \vee 0
 }$.
\end{itemize}
\end{definition}
We now state the analogous of Theorem \ref{mainOV} in the framework
of graph--cycles. The proof is given, for instance, in \cite{Ca}.


\begin{theorem}
\label{mainFW}
(Propositions 4.19, 4.20 and 5.1 \cite{Ca})
Let $A$ be a graph cycle.
For any $\gep > 0$, for any $x,x' \in A$, as $\gb \to \infty$,
one has the asymptotic:
\begin{equation}\label{iptFW}
\bb{P}_{x} \big[\exp\{\gb(H_\rr{e}(A) + \gep) \}
                >
                \tau_{\partial A}
                >
                \exp\{\gb(H_\rr{e}(A) - \gep) \} \big] = 1 - o(1)
\end{equation}
and
\begin{equation}\label{iiptFW}
\bb{P}_{x}\big[\tau_{x'} < \tau_{\partial A},\,
               \tau_{x'} < \exp\{\beta( H_\rr{m}(A) + \gep ) \} \big]
     = 1 - o(1).
\end{equation}
\end{theorem}
We note that this result strongly suggests the equalities
$H_\rr{m}(A) = \Gamma(A)$ and
$H_\rr{e}(A) = \resist(A)$. In the next section we shall prove that
this fact is indeed true.


\subsection{An example.}\label{Esemp}
\par
\begin{figure}[ht!]
  \begin{picture}(200,250)(-80,0)
  \scalebox{0.42}{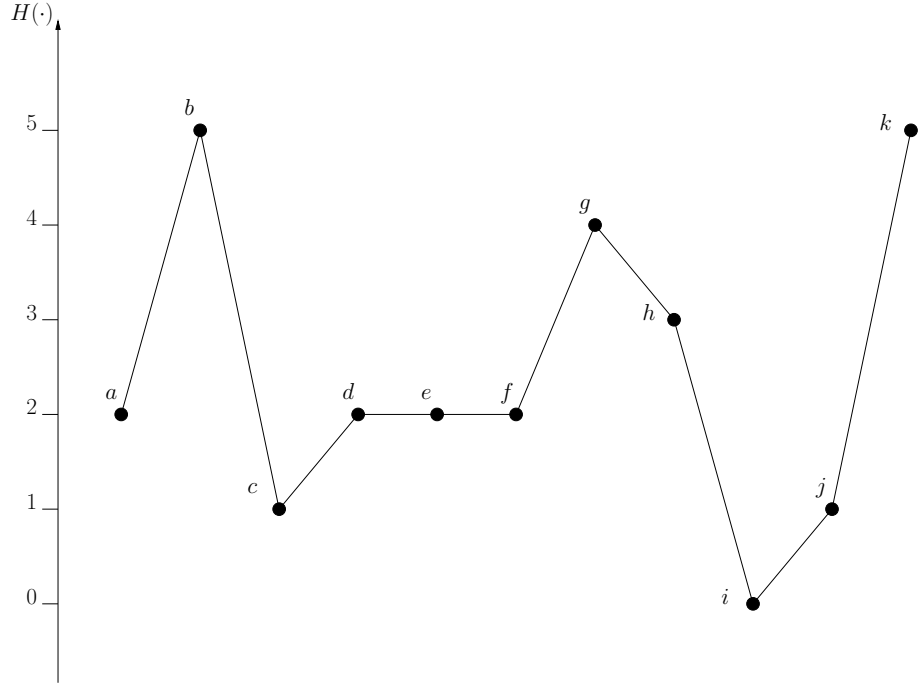}
  \end{picture}
\caption{The energy landscape on $S = \{ a,b,c,d,e,f,g,h,i,j,k\} $}
\label{Fig1}
\end{figure}

\noindent
In this part, we run the algorithm described above in a simple case. We consider the state space $S = \{ a,b,c,d,e,f,g,h,i,j,k\}$ with connectivities and energy landscape described in Figure \ref{Fig1}. For example, $H(a) = H(d) = 2,H(i)=0, q(d,e) > 0, q(e,f) > 0$ and $q(e,g) = 0$.
 
  \textit{Iteration 1.\/} To start our construction, we first recall that 
   \begin{equation*}
    E^{0} = \left\{ \{a\},\{b\},\{c\},\{d\},\{e\},\{f\},\{g\},\{h\},\{i\},\{j\},\{k\} \right\}.
   \end{equation*}

  For two singletons which are not connected (say for example $a$ and $c$), we then have by definition the equality $V^{0}(\{a\},\{c\}) = \infty.$ Else, it is easy to see that $V^{0}(\cdot,\cdot)$ is equal to zero for connected singleton except in the following cases:
   \begin{displaymath}
   V^{0}(\{a\},\{b\}) = 3, V^{0}(\{c\},\{b\}) = 4,V^{0}(\{c\},\{d\}) = 1, V^{0}(\{f\},\{g\}) = 2,  
   \end{displaymath}
   \begin{equation*}
  V^{0}(\{h\},\{g\}) = 1, V^{0}(\{i\},\{j\}) = 1,V^{0}(\{i\},\{h\}) = 3,  
   \;\textrm{ and }\; V^{0}(\{j\},\{k\}) = 4.
   \end{equation*}
    
\par\noindent
    One can then compute the following quantities using \eqref{He}
    \begin{equation*}
    H_\rr{e}^{0}(\{a\})=3,\; H_\rr{e}^{0}(\{c\})=1,\; H_\rr{e}^{0}(\{i\})=1, 
    \end{equation*}
 and 
\begin{displaymath}
H_\rr{e}^{0}(\{b\})\!=\!   
H_\rr{e}^{0}(\{d\})\!=\!  
H_\rr{e}^{0}(\{e\})\!=\!
H_\rr{e}^{0}(\{f\})\!=\!  
H_\rr{e}^{0}(\{g\})\!=\! 
H_\rr{e}^{0}(\{h\})\!=\!  
H_\rr{e}^{0}(\{j\})\!=\!
H_\rr{e}^{0}(\{k\})\!=\!
0.
\end{displaymath}

By making use of equation \eqref{drop} and of the remark after it, 
we note that, for every $x \in S$ and for all $k \geq 0$, 
$ H_\rr{e}^{k}(\{x\}) = H_\rr{e}^{0}(\{x\}) = H_\rr{e}(\{x\})$. 
In other words the remark above gives the value of 
the functions $H_\rr{e}$ and 
$H_\rr{e}^{k}$, for $k \geq 1$, computed at any singleton, namely
$H_\rr{e}(\{x\})$ and 
$H_\rr{e}^{k}(\{x\})$ for $k \geq 1$ and $x\in S$. 
A similar observation will apply in the next steps of the algorithm. 

    Using \eqref{defVs}, we now compute $V^{0}_{\star}(\cdot,\cdot)$, which is infinite for 
not connected singletons and zero for all connected singletons 
apart from the following cases:
    \begin{equation}\label{ex1}
  V^{0}_{\star}(\{c\},\{b\}) = 3,  V^{0}_{\star}(\{f\},\{g\}) = 2, V^{0}_{\star}(\{h\},\{g\}) = 1,
       V^{0}_{\star}(\{i\},\{h\}) = 2, V^{0}_{\star}(\{j\},\{k\}) = 4.
     \end{equation}
Applying the third step of the construction of Section \ref{reccons}, we then deduce that         
  \begin{equation*}
   D^{1} = \left\{ \{a,b\},\{c,d,e,f\},\{g\},\{h\},\{i,j\},\{k\} \right\}.
   \end{equation*} 
For example, to see that $\{a,b\}$ is a subset of $D^{1}$, one notes that $\{a\} \mathcal{R}_{0} \{b\}$ since $V^{0}_{\star}(\{a\},\{b\}) = V^{0}_{\star}(\{b\},\{a\}) = 0$ and $\{b\}$ is not in relation for $\mathcal{R}_{0}$ with $\{c\}$ since  $V^{0}_{\star}(\{c\},\{b\}) = 3 \neq 0$. 
      
  The set $\{a,b\}$ is not a minimal element for  ${\geq}^{1}$ because $\{b\} \xrightarrow{0} \{c\}$, which implies that $\{a,b\} {\geq}^{1} \{c,d,e,f\}$. On the other hand, one checks readily that both $\{c,d,e,f\}$ and $\{i,j\}$ are minimal elements for ${\geq}^{1}$, and we deduce:    
 \begin{equation}\label{Ds}
   D_{\star}^{1} = \{ \{c,d,e,f\},\{i,j\} \}.
   \end{equation}  
   Combining \eqref{En} and \eqref{Ds}, we obtain:
   \begin{equation}\label{Ep}
   E^{1} = \left\{ \{a\},\{b\},\{c,d,e,f\},\{g\},\{h\},\{i,j\},\{k\} \right\}.
   \end{equation}
  \par\noindent  
    We then get immediately by \eqref{Hm} : $H_\rr{m}^{1}(\{x\})= H_\rr{e}(\{x\}) = 0 
    \;\textrm{ for }\; x \in \{b,g,h,k\}$ and 
     \begin{equation*}
      H_\rr{m}^{1}(\{a\})=3,
     H_\rr{m}^{1}(\{c,d,e,f\})=1,
    H_\rr{m}^{1}(\{i,j\})=1.
     \end{equation*}
      
    We can now compute $V^{1}$ making use of \eqref{defVn}, \eqref{Ep} and \eqref{ex1}. Once again, $V^{1}(A,B)$ is infinite as soon as $A$ and $B$ are disconnected, and $V^{1}(A,B) = V^{0}(A,B)$ as soon as $A \in E^{0},B \in E^{0}$. 
For the remaining cases, the cycle pairs with not zero 
$V^1$ are:
  \begin{equation*}
 V^{1}(\{c,d,e,f\},\{b\}) = 4,  V^{1}(\{c,d,e,f\},\{g\}) = 3, V^{1}(\{i,j\},\{h\})= 3, V^{1}(\{i,j\},\{k\})= 5,  
  \end{equation*}
    \begin{equation*}
  V^{1}(\{b\},\{c,d,e,f\}) = V^{1}(\{g\},\{c,d,e,f\}) = V^{1}(\{h\},\{i,j\})= V^{1}(\{k\},\{i,j\})= 0.
  \end{equation*}
  
   \textit{Iteration 2.\/} For the second iteration of the algorithm, one proceeds as in Iteration 1 replacing the state space $E^{0}$ equipped with the cost function $V^{0}$ by the space $E^{1}$ equipped with the cost function $V^{1}$. Applying again \eqref{He} and the remark below \eqref{drop}, we compute:
     \begin{equation*}
     H_\rr{e}(\{i,j\}) = H_\rr{e}^{1}(\{i,j\}) = 3, H_\rr{e}(\{c,d,e,f\})= 3.
     \end{equation*}
    Once again, $V^{1}_{\star}(A,B)$ is infinite as soon as $A$ and $B$ are not connected, identically null otherwise except in the following cases:
     \begin{equation*}
     V^{1}_{\star}(\{c,d,e,f\},\{b\}) = 1, V^{1}_{\star}(\{i,j\},\{k\}) = 2, V^{1}_{\star}(\{h\},\{g\}) = 1, 
     \end{equation*}
      and hence we deduce
    \begin{equation*}
   D^{2} = \left\{ \{a,b\},\{c,d,e,f,g\},\{h,i,j\},\{k\} \right\}.
   \end{equation*} 
    
    Since neither $\{a,b\}$ nor $\{c,d,e,f,g\}$ are minimal elements for the order relation ${\geq}^{2}$, we obtain
   \begin{equation*}
   E^{2} = \left\{ \{a\},\{b\},\{c,d,e,f\},\{g\},\{h,i,j\},\{k\} \right\}.
   \end{equation*}
    Finally we get
  \begin{equation*}
  H_\rr{m}(\{h,i,j\})=3. 
  \end{equation*}  
   
   \textit{Iteration 3.\/} Similarly, for the third iteration of the algorithm, one shows that:
    \begin{equation*}
    H_\rr{e}(\{h,i,j\}) = 4,
    \end{equation*}
    \begin{equation*}
   D^{3} = \left\{ \{a,b\},\{c,d,e,f,g,h,i,j\},\{k\} \right\},
   \end{equation*} 
   \begin{equation*}
   E^{3} = \left\{ \{a\},\{b\},\{c,d,e,f,g,h,i,j\},\{k\} \right\},
   \end{equation*}
  \begin{equation*}
  H_\rr{m}(\{c,d,e,f,g,h,i,j\})=4.
  \end{equation*}
   
    \textit{Iteration 4.\/} After the fourth iteration, the procedure is complete (that is $n_{S} = 4$), and we get: 
     \begin{equation*}
     H_\rr{e}(\{c,d,e,f,g,h,i,j\}) = 5,
     \end{equation*}
   \begin{equation*}
      D^{4} = E^{4} = S,
    \end{equation*}  
    \begin{equation*}
     H_\rr{m}(S) = 5, H_\rr{e}( S) = \infty. 
     \end{equation*} 
     
     We finally deduce the set of graph-cycles of $S$:
    \begin{align*}
     \mathcal{C}(S)  =  \{ \{a\}, \{b\}, \{c\}, \{d\}, \{e\}, \{f\}, \{g\}, \{h\}, \{i\},  \{j\}, \{k\}, \\
      \{c,d,e,f\},\{i,j\}, \{h,i,j\}, \{c,d,e,f,g,h,i,j\}, S \}.
     \end{align*}

     We remark on this toy example that the set of graph cycles coincides with the set of path cycles. We show in the next section that this is indeed always the case.

\section{Equivalence of definitions}
\label{s:equiv}
\par\noindent
In this section we prove that path and graph--cycle decompositions
are equivalent. As a byproduct of the proof we shall also
give a physical interpretation, in terms of differences of energy,
of the quantities $H_\rr{e}$ and
$H_\rr{m}$ defined in the framework of the graph theory.

\begin{theorem}
\label{mainT}
With the notations introduced above:
\begin{enumerate}
\item\label{i:mainT11}
any graph--cycle $A\in\mathcal{C}$ is a path--cycle. Furthermore the following equality holds: 
\begin{equation}
\label{main1}
H_\rr{e}(A)
 = \left( \min_{\partial A} H - \min_{A} H \right ) \vee 0.
  \end{equation}
  One also has 
  \begin{equation}
H_\rr{m}(A) =  \max_{A} H - \min_{A} H  
\end{equation}
 when $|A| > 1$ and $H_\rr{m}(A) = H_\rr{e}(A)$ when $|A| = 1$.
\item\label{i:mainT12}
Any path--cycle is a graph--cycle.
\end{enumerate}
\end{theorem}
 
In particular, starting from anywhere inside a cycle,
the exit time of the cycle on an exponential scale is essentially given by
the depth of the cycle. 
 Note that as soon as $A \in \mathcal{C}$ is a non trivial path cycle, then $H_\rr{e}(A) = \min_{\partial A} H - \min_{A} H$. 
More precisely, we have the equivalence
between Theorems~\ref{mainOV} and \ref{mainFW}.

In order to proof item~\ref{i:mainT11} of the theorem we shall proceed
recursively. In particular we shall introduce
the condition $\mathcal{H}_{n}$ and prove that it holds for any $n \geq 0$.

%

\medskip
\par\noindent
\textit{Condition}~$\mathcal{H}_{n}$. For any $A\in E^{n}$
the following properties hold true:
\begin{displaymath}
\left\lbrace \begin{array}{rl}
       (I)_{n}& A \hspace{2 pt} \text{ is a path cycle;} \hspace{2 pt} \\
       (II)_{n}& \text{if } \hspace{2 pt} A \hspace{2 pt} \text{ (where }  |A| > 1) \text{ and} \hspace{2 pt}  \{a\}     \text{ are connected and belong to } \\
  & E^{n}, \text{ then} \hspace{2 pt}
    V^{n}(A,\{a\}) = H(a) - \min_{A} H \hspace{2 pt}
    \text{ and} \hspace{2 pt} V^{n}(\{a\},A) = 0;  \\
    (III)_{n}& H_\rr{e}^{n}(A) = \left( \min_{\partial A} H - \min_{A} H \right ) \vee 0; \\
       (IV)_{n} & H_\rr{m}^{n}(A) = \max_{A} H - \min_{A} H \hspace{2 pt}
    \text{ when} \hspace{2 pt} |A| > 1.
       \end{array}
\right.
\end{displaymath}

Before turning to the proof of the theorem,
we make a few remarks which hold provided $\mathcal{H}_{k}$ is true for
any $k=0,1,\dots,n-1$,
and that will be used throughout the proof.

R1.\ One important point in the proof will
be to compute $V^{n}(A,B)$ for any element
$A,B \in E^{n}$. Note that if both $|A| > 1$ and $|B| > 1$, then it follows from Lemma \ref{path} that $A$ and $B$ are not connected, and from \eqref{defVn} follows that $V^{n}(A,B) = \infty$. On the other hand, if two singletons $\{a\}$ and $\{b\}$ belonging to $E^{n}$ are connected, then it follows immediately from the definition of $V^{n}$ that $V^{n}(\{a\},\{b\}) = (H(b)-H(a))^{+}$. The other cases are covered by assumption $(II)_{n}$.

R2.\ We state that
      \begin{equation}\label{Dnst}
       E^{n+1} \setminus E^{n} = D_{^\star}^{n+1}.
      \end{equation}
     Indeed, given the definition \eqref{En} of $ E^{n+1}$, for \eqref{Dnst} to hold, we have to show that an element $B \in D_{^\star}^{n+1}$ cannot belong to $E^{n}$. Assume by contradiction that $B \in E^{n}$; there exists at least one element $B' \in E^{n}$ connected to $B$ such that $B \xrightarrow{n} B'$ (and hence $B \geq^{n} B'$). Since $B \in D_{^\star}^{n+1}$, it is a minimal element for $\geq^{n}$, and hence one gets $B' \geq^{n} B$. This implies that both $B$ and $B'$ are in the same equivalence class for $\mathcal{R}_{n}$, which is contradictory in view of the construction of $D_{^\star}^{n+1}$.

     This remark is very useful, indeed,
to prove the implication $\mathcal{H}_{n} \Rightarrow \mathcal{H}_{n+1}$, we   show the four properties defining $\mathcal{H}_{n+1}$ restricting
the analysis to the case $A \in D_{^\star}^{n+1}$.

R3.\ We state an important decomposition of any $A\in E^{n+1}\setminus E^n$
which may be viewed as the
analogous of \cite[Proposition~6.19]{OV} in the graph--cycle context.

       We write $A$ as the disjoint union
      \begin{equation}\label{dec}
       A = \bigsqcup_{j=1}^{j_{A}} A_{j}
      \end{equation}
      where $j_{A} \geq 1$ and,
for all $j \leq j_{A}, A_{j} \in E^{n}$. This decomposition is in fact the decomposition $\mathcal{M}_{\star}(A)$ which we rewrite in a more tractable way.

         Of course $A$ not belonging to $E^{n}$ implies in fact that $j_{A} \geq 2$. In the remaining of the proof, we will refer to the elements appearing in the decomposition of the right hand side of \eqref{dec} as \textit{subcycles} of $A$.

        For $j \in [1,j_{A}]$, whenever there exists $a \in S$ such that $A_{j} = \{a\}$ and $b \in S$ a neighbor of $a$ such that $H(b) \leq H(a)$,
we say that the singleton $A_{j}$ is \textit{a good singleton of $A$}. In words, a good singleton $\{a\}$ of $A$  is a trivial path cycle because $a$ is not a local minimum of $H(\cdot)$. We denote by $j_{A}^{s}$ the number of subcycles of $A$ which are good singletons of $A$. Note that whereas the notion of trivial path cycles depends only on the energy landscape, the notion of good singleton of $A$ strongly depends on the cycle decomposition of $A$.

R4.\  $(I)_{n}$ implies that, as soon as $A_{j} \xrightarrow{n} \{a\}$, where $|A_{j}| > 1$ and $\{a\}$ is a subcycle of $A$ such that $a \in \partial A_{j}$, then necessarily $\{a\}$ is good.






R5.\ For singletons which are not good,
say that $A_{j} = \{a\}$ where $j > j_{A}^{s}$,
we can remark that all the subcycles of $A$ which are connected to $\{a\}$ are good singletons and have the same energy.
  Indeed, all these subcycles are path cycles by $(I)_{n}$. Since a non--trivial path cycle cannot be connected to a local minimum for $H(\cdot)$ in view of \eqref{defC}, it follows that all these subcycles are good singletons. To see that they have the same energy, consider any  $b \in \partial \{a \} \cap A$ satisfying $\{a\} \xrightarrow{n} \{b\}$. By the definition of $\xrightarrow{n}$, we then get that
$H(b) = H(a) + H_\rr{e}(\{a\}).$

 R6 Combining R4, R5 and Lemma \ref{path}, we note the important fact that $j_{A}^{s} \geq 1$. Up to reordering, we assume from now on that for $1 \leq j \leq j_{A}^{s}$, $A_{j}$ is a good singleton.

R7.\ We then show the following equalities:

        \begin{equation}
        \label{singls}
         H\left ( \bigcup_{1 \leq j \leq j_{A}^{s}}A_{j}\right ) = \max_{A} H
         \end{equation}
         and
        \begin{equation}\label{singlsb}
         H\left (\partial A_{j} \cap A, j \in [j_{A}^{s}+1,j_{A}] \right ) = \max_{A} H .
        \end{equation}

          Of course, equations \eqref{singls} and \eqref{singlsb} implicitly state that $H(A_{j})$ does not depend on $j \in [1,j_{A}^{s}]$ and that $H\left (\partial A_{j} \cap A \right)$ does not depend on $j \in [j_{A}^{s}+1,j_{A}]$.

           To prove \eqref{singls} and \eqref{singlsb}, we distinguish the cases $j_{A}^{s} = 1$ and $j_{A}^{s} \geq 2$.

           If $j_{A}^{s} = 1$, then for every $j \in [2,j_{A}]$, it follows from Lemma \ref{s:path} and remark R4 that $A_{j}$ is connected to $A_{1}$ and $A_{j} \xrightarrow{n} A_{1}$ for $j \in [2,j_{A}]$. As a consequence, considering $(II)_{n}$, we get both \eqref{singls} and \eqref{singlsb}.

         Assume now that  $j_{A}^{s} \geq 2$ and consider two elements $A_{i}$ and $A_{j}$ (where $1 \leq i,j \leq j_{A}^{s}$ and $i \neq j$). We consider a path $\omega = (\omega_{1}, \dots, \omega_{n})$ of connected subcycles of $A$ such that for $\omega_{1} = A_{i}$, $\omega_{n} = A_{j}$ and for any $k \in [1,n-1], \omega_{k} \xrightarrow{n} \omega_{k+1}$.

           Let us consider $k_{0} \in [1,n-2]$ such that $\omega_{k_{0}} = A_{l}$ where $l \leq j_{A}^{s}$, or in words $\omega_{k_{0}}$ is a good singleton of $A$. We claim that either $\omega_{k_{0}+1}$ or $\omega_{k_{0}+2}$ is a good singleton of $A$, and moreover that in both cases the energy of this good singleton is equal to the energy of  $\omega_{k_{0}}$.

            Indeed, if $\omega_{k_{0}+1}$ is not a good singleton, then either $|\omega_{k_{0}+1}| > 1$ or $\omega_{k_{0}+1}$ is a local minimum for $H$. When $|\omega_{k_{0}+1}| > 1$, it follows  from remark R4 that necessarily  $\omega_{k_{0}+2}$ is a good singleton. Also, it follows from $(II)_{n}$ that $H(\omega_{k_{0}}) = H(\omega_{k_{0}+2})$. On the other hand, in the case where $\omega_{k_{0}+1}$ is a local minimum for $H$, we saw in remark R5
that necessarily $\omega_{k_{0}+2}$ is a singleton and that its energy is the same as the one of $\omega_{k_{0}}$.

     It follows easily from these considerations that $H(A_{i}) = H(A_{j})$, and hence $H(A_{j})$ does not depend on $j \in [1,j_{A}^{s}]$.

      The equalities \eqref{singls} and \eqref{singlsb} then follow readily from Lemma \ref{path}, from $(II)_{n}$ and $(I)_{n}$.

\par\noindent
\textit{Proof of Theorem \ref{mainT}.\/}
Proof of item~\ref{i:mainT11}: we shall prove recursively that
Condition $\mathcal{H}_n$ holds for any $n\ge0$, which
implies item~\ref{i:mainT11}. 
  Since singletons are path--cycles, $\mathcal{H}_{0}$ is trivially true. Note in particular that for any $x \in S$, the equality $H_{\rr{m}}(\{x\}) = H_{\rr{e}}(\{x\})$ holds true.

Let us now assume that $\mathcal{H}_{n}$ holds for a given $n \geq 0$ and
prove that
$\mathcal{H}_{n+1}$ holds. Considering R2, we restrict ourselves to the case where $A \in D_{\star}^{n+1}$ (and in particular $|A| > 1$); in the remaining of the proof of Theorem \ref{mainT}, $A$ will be such an element.

         $(IV)_{n+1}$  is a consequence of \eqref{singls}, $(III)_{n}$ and $(II)_{n}$. Indeed, since $|A| > 1$, combining the definition \eqref{Hm} of $H_\rr{m}$ and the induction hypothesis, we get the identity
      \begin{equation}\label{IVn}
       H_\rr{m}^{n+1}(A) = \max \{ \left (\min_{\partial A'} H - \min_{A'} H\right ) \vee 0, A' \in E^{n}, A' \subset A \}.
      \end{equation}

       In the case where $j_{A}^{s} = j_{A}$, then \eqref{IVn} immediately gives that $H_\rr{m}(A)=0$, which implies $(IV)_{n+1}$ using \eqref{singls}.

        Assume now that $j_{A} > j_{A}^{s}$. We  note that one can restrict the optimization appearing in the right hand side of \eqref{IVn} to  subcycles $A'$ of $A$ verifying $A' = A_{j}$ for $j \in [j_{A}^{s}+1;j_{A}]$; indeed, it follows from $(I)_{n}$ (in the case where $|A_{j}| > 1$) and from remark R5
(when $|A_{j}| = 1$) that the quantity $\min_{\partial A'} H - \min_{A'} H$ is positive for such cycles and it is identically null as soon as $A'$ is a good singleton.

    For such an $A'$, the equality \eqref{singlsb} implies that the quantity $\min_{\partial A'} H$ does not depend on $A'$ and is equal to $ \max_{A} H$. Hence one gets
       \begin{equation*}
       H_\rr{m}^{n+1}(A) = \max_{A} H - \min_{  j \in [j_{A}^{s}+1,j_{A}] } \min_{A_{j}} H.
      \end{equation*}

      Since $E^{n}$ is a partition of $S$ and since $\max_{A} H $ is reached on $A \setminus \bigcup_{1 \leq j \leq j_{A}^{s}} A_{j}$ as shown in \eqref{singls}, one gets  $\min_{   j \in [j_{A}^{s}+1,j_{A}] } \min_{A_{j}} H = \min_{A} H$, and hence the equality $(IV)_{n+1}$ holds in this case as well.

      $(III)_{n+1}$ is a direct consequence of $(II)_{n}$ and $(I)_{n}$. By Lemma \ref{path}, we get that all the neighbors of $A$ in $E^{n}$ are singletons, and hence we have:
    \begin{align*}
H_\rr{e}^{n+1}(A) & = \inf_{\{a\} \in E^{n}, a \in \partial A} V^{n}(A,\{a\})
         = \inf_{\{a\} \in E^{n}, a \in \partial A} H(a) - \min_{A} H \\
          &  = \min_{\partial A} H - \min_{A} H
          \end{align*}
           where in the second equality we made use of $(II)_{n}$.

            $(I)_{n+1}$ is a consequence of $(III)_{n+1}$ and $(IV)_{n+1}$. Indeed, we have to show that  $ \max_{A} H < H(F( \partial A))$. Note that this is equivalent to showing that $H_\rr{m}^{n+1}(A) < H_\rr{e}^{n+1}(A)$, which in turn by definition
        of $H_\rr{e}(\cdot)$ is equivalent to the fact that, for every $A' \in E^{n+1}$:
         \begin{equation*}
          \min\{ V_{\star}^{n}(B,B'), B,B' \in E^{n}, B \subset A, B' \subset A' \} > 0.
         \end{equation*}

          By contradiction, assume that there exist $ A' \in E^{n+1}, A' \neq A,  B,B' \in E^{n}, B \subset A, B' \subset A'$ such that $V_{\star}^{n}(B,B') = 0$ (and of course $A \neq A'$); this implies that $A \geq^{n+1} A'$, which contradicts the minimality of $A$ for the order relation $\geq^{n+1}$ (given the construction of $A$, this minimality is necessarily strict).

       $(II)_{n+1}$ is a consequence of $(IV)_{n+1},(III)_{n+1}$ and $(II)_{n}$. Indeed, for $a \in \partial A$ such that $\{a\} \in E^{n+1}$, using the definition \eqref{defVn} of $V_{n}$, of $(IV)_{n+1}$ and of $(I)_{n}$, we get that:
        \begin{equation}\label{recV}
            V^{n+1}(A,\{a\})  = \max_{A} H - \min_{A} H + \min\{ V_{*}^{n}(B,\{a\}), B \in E^{n}, B \subset A \}.
               \end{equation}
                Note that we restricted the set on which we minimize  $V_{*}^{n}$  in the right hand side of  \eqref{recV}  since in all other cases this quantity is infinite.

          Now we distinguish two cases.

           In the case where the $\min$ in the right hand side of  \eqref{recV} is attained on a subcycle $A_{j_{0}}$ of $A$ such that $ j_{0} > j_{A}$, making use of $(II)_{n}$ and $(III)_{n}$, we get:
          \begin{equation}\label{cascy}
      \min\{ V_{*}^{n}(B,\{a\}), B \in E^{n}, B \subset A \} = H(a) - \min_{A_{j_{0}}} H - (\min_{\partial A_{j_{0}}} H - \min_{A_{j_{0}}} H) = H(a) - \min_{\partial A_{j_{0}}} H.
               \end{equation}
         Then we already noted that $\min_{\partial A_{j_{0}}} H = \max_{A} H$ holds because $ j_{0} > j_{A}$ and of equality \eqref{singlsb}.

          When the $\min$ in the right hand side of  \eqref{recV} is attained on a good   singleton $A_{j_{1}}$ with $j_{1} \in [1,j_{A}^{s}]$, one gets that
          \begin{equation}
          \label{cassi}
       \min\{ V_{*}^{n}(B,\{a\}), B \in E^{n}, B \subset A \} = (H(a)-H(A_{j_{1}}))^{+} = H(a)-  \max_{A} H.
          \end{equation}

       Equation \eqref{cassi} holds because $A$ satisfies \eqref{defC}, and hence $(H(a)-H(A_{j_{1}}))^{+} = H(a)-H(A_{j_{1}})$, and on the other hand $H(A_{j_{1}}) = \max_{A} H$ from \eqref{singls}.

        Finally, we have that 
$V^{n+1}(\{a\},A) = H_{m}(\{a\}) + \min\{ V^{n}(\{a\},A_{j}) - H_{e}(\{a\}), j \in [1,j_{A}] \} = \min\{ V^{n}(\{a\},A_{j}), j \in [1,j_{A}] \}$. When this minimum is attained on a non--trivial subcycle of $A$, one gets $V^{n+1}(\{a\},A) = 0$ by $(II)_{n}$, and when it is attained on a singleton $\{b_{0}\}$ involved in the decomposition \eqref{dec}, using $(I)_{n+1}$, one gets that $V^{n+1}(\{a\},A) =((H(b_{0}) -H(a))^{+} = 0 $, thus in any case
         \begin{equation}
           \label{cassicy}
            V^{n+1}(\{a\},A) = 0.
           \end{equation}

        Combining \eqref{recV}, \eqref{cassi}, \eqref{cascy} and \eqref{cassicy}, we get $(II)_{n+1}$, and hence the recursion is completed. As noted
above this completes the proof of item~\ref{i:mainT11}.

     Now we prove item~\ref{i:mainT12} of the theorem. Let us
     consider $A$ a connected subset of $S$ which is not a singleton, and such that $\max_{A} H < \min_{\partial A} H$.

    We introduce $k_{0} := \inf\{ k \geq 1, \exists A' \in E^{k}, A \subset A' \}$ and we denote by $A'$ the element of $E^{k_{0}}$ such that $A \subset A'$. Since the recursion of section \ref{s:FWC} is not stationary, $k_{0}$ is well defined.
 To show Theorem \ref{mainT}, it is enough to show that the reciprocal inclusion holds.
     We introduce the decomposition of $A'$ into disjoint elements of $E^{k_{0}-1}$ as we did in \eqref{singls}:
      \begin{equation}\label{decb}
       A' = \bigsqcup_{j=1}^{j_{A'}} A'_{j},
      \end{equation}
       and we define $j_{A'}^{s}$ for $A'$ as we defined $j_{A}^{s}$ for $A$ in R3.

      The inclusion $A' \subset A$ will hold
     as soon as we show that there exists an element $a \in S$ such that $\{a \}$ is a good singleton for $A'$ and $a \in A$; indeed, we already proved in item \ref{i:mainT11} that $A'$ is a path cycle, hence combining Lemma \ref{lvlset} and equality \eqref{singls}, we remark that $A' = U_{\leq a}$. From Lemma \ref{interincl}, we thus get $A' \subset A$.

    By the minimality of $k_{0}$, $A$ is not contained in any of the $(A'_{j})_{j \leq j_{A'}}$, and thus there exist
    $j_{1}$ and $j_{2}$ in $[1;j_{A'}]$ such that both $A \cap A'_{j_{1}} \neq \emptyset$ and  $A \cap A'_{j_{2}} \neq \emptyset$. $A$ being connected, there exists a path $\omega = (\omega_{1}, \ldots, \omega_{n})$
     contained in $A$ joining $A \cap A'_{j_{1}}$ to $A \cap A'_{j_{2}}$ (say $\omega_{1} \in A_{j_{1}}$ and $\omega_{n} \in A_{j_{2}}$); we then consider three cases:
 \begin{enumerate}
      \item either $A'_{j_{1}}$ or $A'_{j_{2}}$ is a good singleton, and then we are done.
      \item either $A_{j_{1}}$ or  $A_{j_{2}}$ is not a singleton, and then we assume by symmetry that $|A_{j_{1}}| > 1$. Let us consider $l = \inf_{j \geq 1} \{ \omega_{j} \notin A_{j_{1}} \}$; by Lemma \ref{path}, $\{ \omega_{l} \}$ is necessarily a subcycle of $A'$, and by R4, $\{ \omega_{l} \}$ is a good singleton.
       \item in the case where both $A_{j_{1}}$ and  $A_{j_{2}}$ are not good singletons, then it follows from R5 that they cannot be connected and hence $n \geq 3$. It follows from R5 as well that $\{ \omega_{2} \}$ is necessarily a good singleton.
      \end{enumerate}

    \qed



\end{document}